\documentclass[11pt]{amsart}
\usepackage{amssymb}
\usepackage{stmaryrd}
\usepackage[all]{xy}
\usepackage{epsf}

\newtheorem{theorem}{Theorem}
\newtheorem{proposition}{Proposition}
\newtheorem{lemma}{Lemma}
\newtheorem{corollary}{Corollary}
\newtheorem{conjecture}{Conjecture}

\newtheorem{definition}{Definition}

\theoremstyle{definition}

\theoremstyle{remark}
\newtheorem{remark}{Remark}

\newcommand{\rmap}{\longrightarrow}
\DeclareMathOperator{\Der}{Der}         
\DeclareMathOperator{\Def}{def}         %
\newcommand{\de}[1]{{#1}_{\Def}}        
\newcommand{\D}{\ensuremath{\mathcal{D}}}
\newcommand{\tto}{\rightrightarrows}
\newcommand{\al}{\alpha}                
\newcommand{\be}{\beta}                 
\DeclareMathOperator{\ad}{ad}         

\newcommand{\comment}[1]{}

\begin{document}

\title[Deformations of Lie brackets: cohomological aspects]{Deformations of Lie brackets: cohomological aspects}

\author{Marius Crainic}
\address{Depart. of Math., Utrecht University, 3508 TA Utrecht, 
The Netherlands}
\email{crainic@math.uu.nl}
\author{Ieke Moerdijk}
\address{Depart. of Math., Utrecht University, 3508 TA Utrecht, 
The Netherlands}
\email{moerdijk@math.uu.nl}
\thanks{The first author was partially supported by a  KNAW Fellowship (Utrecht) and a Miller
Fellowship (Berkeley). }


\begin{abstract}
We introduce a new cohomology for Lie algebroids, and prove that it
provides a differential graded Lie algebra which ``controls'' 
deformations of the structure bracket of the algebroid. 
\end{abstract}
\maketitle

\section{Introduction}             %
\label{Deformation: Outline}      %

The aim of this paper is to find the differential graded Lie algebra and its cohomology theory controlling
deformations of a large class of geometric structures, known as Lie algebroids. This problem is particularly 
difficult since the notion of adjoint representation is not available for these structures. 

We recall that a Lie algebroid over a manifold $M$ is a vector bundle
$\pi:A\to M$ together with a Lie bracket $[~,~]$ on the space of
sections $\Gamma(A)$ and a bundle map $\rho :A\to TM$, called the anchor,
satisfying the Leibniz identity:
\[ [\al,f\be]=f[\al,\be]+L_{\rho\al}(f)\be, \qquad f\in
C^\infty(M),\al,\be\in\Gamma(A).\]

The notion of Lie algebroid goes back in this form to Pradines in 1967, but, in local coordinates, it already 
appeared in E. Cartan's work in 1904, and the analogous algebraic notion was already studied by Rinehart in 1963.

Lie algebroids are to be thought of as infinite-dimensional Lie
algebras of ``geometric type'', or as generalized tangent bundles.
Indeed, the simplest examples are (finite-dimensional) Lie algebras and
tangent bundles of manifolds, and there are many natural geometric
examples coming from foliations, Poisson manifolds, infinitesimal
actions of Lie algebras on manifolds, and other contexts.
The present paper is thus part of a larger programme, the goal of which
is to give a unified approach to geometric structures controlled by
``Lie brackets'', and to make explicit the analogies and interplay between
the various fields such as Lie group theory, theory of foliations,
Poisson geometry, etc.

A well-known principle in mathematics, present already in the study of deformations of complex structures (cf. \cite{GM}
and the references therein) and in the work of Nijenhuis-Richardson \cite{NR1}, and
emphasized by P. Deligne, M. Kontsevich, and
others, states that the deformation theory of a specific (type of)
structure is governed by a naturally associated differential graded Lie
algebra (dgla). Moreover, this dgla, or the graded Lie algebra given by
its cohomology, should act on invariants associated to the structure.
In this paper, we will exhibit for any Lie algebroid $A$ the dgla
governing its deformations, the cohomology Lie algebra of which we call
the deformation cohomology of $A$ and denote by $\de{H}^{*}(A)$. We will
relate this cohomology to known cohomologies of $A$, and in particular
prove that it acts on the ``De Rham'' cohomology of $A$. We will prove
that deformations of the bracket of a given algebroid $A$ give rise to
cohomology classes in $\de{H}^{2}(A)$, and we investigate the relation
of these classes to those constructed earlier for the deformation of
foliations \cite{Hei} and of Poisson brackets. We will also prove that, quite surprisingly,
the vanishing of such cohomology classes do imply rigidity (Theorem \ref{rigidity}).

In the extreme case where the Lie algebroid $A$ is simply a finite
dimensional Lie algebra, it is of course well-known that the
deformations of the Lie bracket are controlled by (the cohomology of)
the Chevalley-Eilenberg complex for the adjoint representation of $A$ \cite{NR1, NR2}.
Our theory contains this fact as a special case, and can be interpreted
as a way of defining the cohomology of a general Lie algebroid ``with
coefficients in its adjoint representation'', in spite of the fact that
this adjoint representation does not exist in the general context of
Lie algebroids. (Several efforts have been made to remedy this
situation \cite{coh}, but no satisfactory answer has been found so far.)

Based on the analogy with the adjoint representation and the rigidity
properties of compact Lie groups, we state at the end of this
paper a general ``rigidity conjecture'' for algebroids of compact
type which we expect to have applications to the linearization problem \cite{We0} and other 
rigidity problems.
Moreover, we prove this conjecture in two special cases, namely when the Lie
algebroid is regular and
when the Lie algebroid is defined from an infinitesimal action of a Lie
algebra on a manifold. The first case is closely 
related to Weinstein's linearization theorem in the regular case
\cite{We}- note in particular that his proof uses the same techniques, based on a Bochner-type averaging
and the Van Est isomorphism for groupoids \cite{Cra}. The second case is related 
to the linearization theorem for Poisson manifolds \cite{Conn}. 
We also give an idea of a proof in the general
case, based on a ``lin''-version of the category of smooth manifolds,
which we expect to be of independent interest.

{\bf Acknowledgements:} First drafts of this paper were written while the second author was
visiting the University of Ljubljana (November 2000), and later while
the second author was a Miller Fellow at UC Berkeley (2001). The main
results of this paper were first presented during the Foliation Theory
Program at the ESI in Vienna (October 2002). We are grateful to 
F. Kamber
and K. Richardson for inviting us to this program, and to
K. Mackenzie, J. Mrcun, D. Roytenberg, A. Weinstein, and especially to R.L.
Fernandes, for helpful discussions related to this paper. After having posted the paper
on the Archive, Y. Kosmann-Schwarzbach kindly pointed out several innacuracies to us, while
J. Grabowski informed us that the definition of the deformation complex also occurs in the recent 
preprint DG/0203112v2.

\section{Multiderivations and the deformation complex}%
\label{Deformation complex}       %

The {\it deformation complex} of a Lie algebroid $A$ is defined as the complex
$(\de{C}^{*}(A),\delta)$ in which the $n$-cochains $D\in \de{C}^{n}(A)$ are multilinear
antisymmetric maps
\[ D: \underbrace{\Gamma(A)\times \ldots \times \Gamma(A)}_{n\ \text{times}} \rmap \Gamma(A)\]
which are multiderivations, and the coboundary is given by
\begin{multline}
\label{differential}
\delta(D)(\alpha_0,\dots,\alpha_{n})=
\sum_{i} (-1)^i [\alpha_i,D(\alpha_0,\dots,\hat{\alpha_i},\dots,\alpha_{n})]+\\
+\sum_{i<j} (-1)^{i+j}
D([\alpha_i,\alpha_j],\alpha_0,\dots,\hat{\alpha_i},\dots,\hat{\alpha_j},\dots,\alpha_{n}).
\end{multline} 
We denote by $\de{H}^{*}(A)$ the resulting cohomology.
In this section we explain this definition and we investigate the structure
present on $(\de{C}^{*}(A),\delta)$.

\subsection{Multiderivations}
\label{Multiderivations}
Let $E\to M$ be a vector bundle, and denote by $r$ its rank. Recall
that a {\it derivation} on $E$ is any linear operator
$D: \Gamma(E)\to \Gamma(E)$
with the property that there exists a vector field $\sigma_D\in\mathcal{X}(M)$,
called the {\it symbol} of $D$, such that
\[ D(fs)= fD(s)+ \sigma_D(f)s, \]
for any section $s\in \Gamma(E)$ and function $f\in C^{\infty}(M)$. 
We denote by $Der(E)$ the space of derivations on $E$.

Assume, for the moment, that the rank of $E$ is $r\ge 2$.  A
{\it multiderivation of degree n} is a skew-symmetric multilinear
map
\[ 
D: 
\underbrace{\Gamma(E)\otimes\cdots\otimes\Gamma(E)}_{n+1\ \text{times}}
\to \Gamma(E) 
\]
which is a derivation in each entry, i.e., there is a map
\[ 
\sigma_{D}: 
\underbrace{\Gamma(E)\otimes\cdots\otimes\Gamma(E)}_{n\ \text{times}}
\to \mathcal{X}(M), 
\] 
called the {\it symbol} of $D$, such that
\[ D(s_0,s_1,\dots, fs_n)= 
fD(s_0,s_1,\dots,s_n)+\sigma_{D}(s_0,\dots,s_{n-1})(f)s_n \] 
for any function $f\in C^{\infty}(M)$ and sections $s_i\in\Gamma(E)$. Notice that this
identity determines $\sigma_{D}$ uniquely.

We will denote by $Der^n(E)$ the space of multiderivations of degree
$n$, $n\geq 0$. We have $Der^0(E)= Der(E)$, and we set $Der^{-1}(E)= \Gamma(E)$. 

\begin{lemma}
\label{deriv} 
For any multiderivation $D\in Der^n(E)$ of degree $n\geq 0$, its symbol $\sigma_D$ is
anti-symmetric and $C^{\infty}(M)$-linear. Moreover, $Der^n(E)$ is
the space of sections of a vector bundle $\D^{n}E\to M$ which fits
into a short exact sequence of vector bundles
\[ 
\xymatrix{
0\ar[r]& \wedge^{n+1}E^{\vee}\otimes E\ar[r]&
\D^{n}E\ar[r]& \wedge^{n}E^{\vee}\otimes TM \ar[r]& 0,}
\] 
where $E^{\vee}$ denotes the dual of $E$. In particular, $Der^n(E)= 0$ for $n\geq \text{rk}(E)+ 1$. 
\end{lemma}

\begin{proof} 
The antisymmetry of $\sigma_D$ follows from that of $D$. From this it
follows also that
\[ 
D(s_0,s_1,\dots , s_{i-1}, fs_i, s_{i+1}, \ldots ,s_n)=fD(s_0,s_1,\dots,s_n)+(-1)^{n-i}\sigma_{D}(s_0,\dots ,\widehat{s_i}, \ldots ,s_n)(f)s_i.
\]
Now compute $D(fs_0, gs_1, \ldots , s_n)$ in two ways: first by 
taking $f$ out followed by taking out $g$, and then the other way. We
obtain 
\begin{multline*}
\left(\sigma_D(fs_0,s_2,\dots,s_n)(g)-f\sigma_D(s_0,s_2,\dots,s_n)(g)\right)s_1+ \\
+\left(\sigma_D(gs_1,s_2,\dots,s_n)(f)-g\sigma_D(s_1,s_2,\dots,s_n)(f)\right)s_0=0.
\end{multline*}
Since $E$ was assumed to be of rank $r\geq 2$, it follows that
\[ \sigma_D(fs_0,s_2,\dots,s_n)- f\sigma_D(s_0,s_2,\dots,s_n)=0,\] 
i.e., $\sigma_D$ is $C^{\infty}(M)$-linear.

Observe now that, if $\nabla$ is a connection on $E$, the operator 
\[ L_{D}(s_{0},\dots,s_{n})= D(s_{0},\dots,s_{n})+ 
(-1)^n \sum_{i}(-1)^{i+1}\nabla_{\sigma_{D}(s_{0},\dots,\hat{s}_{i},\dots,s_{n})}(s_i) 
\] 
is antisymmetric and $C^{\infty}(M)$-multilinear. This shows that any
connection $\nabla$ on $E$ determines an isomorphism of
$C^{\infty}(M)$-modules 
\[ Der^n(E)\cong \Gamma(\wedge^{n+1}E^{\vee}\otimes E)\oplus
\Gamma(\wedge^{n}E^{\vee}\otimes TM),\]
which sends $D$ to $(L_D, \sigma_D)$. This proves the statement about the
existence of a vector bundle $\D^n E$, and it shows at the same time
that the choice of a connection $\nabla$ induces a splitting of the
short exact sequence above.
\end{proof}

A similar definition when the rank of $E$ is $r=1$ would not imply
the $C^{\infty}(M)$-linearity of the symbols. To extend the definition
of $Der^n(E)$ to this case, we must require $C^{\infty}(M)$-linearity 
of its symbols.

\subsection{Multiderivations and brackets}  
There is a close connection between the spaces $Der^n(E)$ and
Lie algebroids. First of all, the obvious bracket $[\cdot, \cdot]$ on
$Der(E)$ extends to  $Der^*(E)$ as follows:  

\begin{proposition} 
\label{brackets}
For $D_1\in Der^{p}(E)$ and $D_2\in Der^q(E)$ we define the Gerstenhaber bracket 
\[ [D_1, D_2]= (-1)^{pq}D_1\circ D_2- D_2\circ D_1 ,\]
where 
\begin{eqnarray}
\label{gert-eq}
 & D_2\circ D_1 (s_0,\dots,s_{p+q})  = & \nonumber \\
 &  \sum_{\tau} (-1)^{|\tau|} 
D_2(D_1(s_{\tau(0)},\dots,s_{\tau(p)}),s_{\tau(p+1)},\dots,s_{\tau(p+q)}))  & 
\end{eqnarray}
and the sum is over all $(p+1, q)$-shuffles. Then $[D_1, D_2]\in Der^{p+q}(E)$, and
the resulting bracket $[\cdot, \cdot]$ makes $Der^*(E)$ into a graded Lie
algebra.
\end{proposition}

\begin{proof}
It is well-known that the bracket above provides a gla structure on the space 
of skew-symmetric multilinear maps $\Gamma(E)\otimes \ldots \otimes \Gamma(E)\rmap \Gamma(E)$ \cite{NR2}. Therefore, it
suffices to show that the space of multiderivations is closed under this bracket, and this follows 
by a careful (but straightforward) calculation. The final conclusion of the computation is that $[D_1, D_2]$ is a multiderivation
with symbol:
\[ \sigma_{[D_1, D_2]}= ((-1)^{pq} \sigma_{D_1}\circ D_2- \sigma_{D_2}\circ D_1)+ [\sigma_{D_1}, \sigma_{D_2}] ,\]
where
$\sigma_{D_1}\circ D_2$ is given by the same formula as (\ref{gert-eq}) above (so $\sigma_{D_1}\circ D_2= 0$ if $p=0$), while
\[ [\sigma_{D_1}, \sigma_{D_2}](s_1, \ldots , s_{p+q})= \sum_{\tau} (-1)^{|\tau|} 
[\sigma_{D_1}(s_{\tau(1)}, \ldots , s_{\tau(p)}), \sigma_{D_2}(s_{\tau(p+1)}, \ldots , s_{\tau(p+q)}]\]
(sum over shuffles again). 
\end{proof}

The formulas in Proposition \ref{brackets} are quite standard and they go back to 
Gerstenhaber \cite{Ger} (the case of algebras) and Nijenhuis and Richardson \cite{NR2, NR3} (the case of Lie algebras). 
Note however that we have chosen the signs differently, so as to match the classical formulas for Lie derivatives
and the De Rham differential.

The gla structure in $Der^*(E)$ allows us to give the following
folklore description of Lie algebroids (going back at least to \cite{Ger}):

\begin{lemma}
\label{rer-alg}
If $A$ is a vector bundle over $M$, then there exists a one-to-one
correspondence between Lie algebroid structures on $A$ and elements
$m\in Der^1(A)$ satisfying $[m,m]=0$.
\end{lemma}

The more familiar form of the definition of a Lie algebroid is
obtained by letting $[\alpha,\beta]= m(\alpha,\beta)$ (the Lie bracket) and
$\rho=\sigma_{m}:A\to TM$ (the anchor).

The vector bundle $\D^0 E$ itself is a Lie algebroid for any vector
bundle $E$. The bracket is the one mentioned above (given by the
commutators of derivations), while the anchor is just taking the
symbol, $\rho(D)= \sigma_{D}$.

Let us recall \cite{Ma} that a representation of a Lie algebroid $A\to M$ is a vector
bundle $E\to M$ together with a flat $A$-connection $\nabla$ on $E$. This means
that $\nabla:\Gamma(A)\otimes \Gamma(E)\to \Gamma(E)$ is a bilinear
map, written $(\alpha,s)\mapsto \nabla_{\alpha}(s)$, which satisfies the 
connection properties
\begin{align*}
\nabla_{f\alpha}(s)&=f\nabla_{\alpha}(s),\\ 
\nabla_{\alpha}(fs)&=f\nabla_{\alpha}(s)+ \rho\alpha(f)s,
\end{align*} 
as well as the flatness condition 
\[\nabla_{[\alpha, \beta]}= [\nabla_{\alpha}, \nabla_{\beta}].\]
We observe the following well-known fact:

\begin{lemma} 
\label{lemma:Lie algebroid}
Given a Lie algebroid $A$ over $M$, there exists a one-to-one
correspondence between representations of $A$ and vector bundles $E$
over $M$ together with a Lie algebroid map $\nabla: A\to \D^0 E$. 
\end{lemma}

The lemma suggests using the notation $\mathfrak{gl}(E)$ for
the Lie algebroid $\mathcal{D}^0E$. 
Of course, $E$ is a representation of 
$\mathfrak{gl}(E)$, with the tautological action 
\[ \nabla_{D}(s)=D(s),\] 
for $D\in D^0(E)$ and $s\in \Gamma(E)$. 
We point out that the Lie algebroid
$\mathfrak{gl}(E)$ is always integrable: one checks easily that the
Lie groupoid $GL(E)\tto M$ for which the arrows
$x\stackrel{g}{\rmap}y$ are the linear isomorphisms $g:E_x\to E_y$, has
Lie algebroid precisely $\mathfrak{gl}(E)$ (see also \cite{Ma}).

\subsection{Cohomology}
\label{coh-subsect}
Recall that, given a representation $E= (E, \nabla)$ of a Lie algebroid $A$,
the {\it De Rham cohomology} of $A$ with coefficients in $E$ \cite{Ma}, denoted $H^{*}(A; E)$, is
defined as the cohomology of the complex $C^{*}(A; E), \delta_{A, E})$, where
$C^p(A; E)= \Gamma(\Lambda^pA^{\vee}\otimes E)$ consists of 
$C^{\infty}(M)$-multilinear antisymmetric maps 
\[ \underbrace{\Gamma(A)\times \ldots \times \Gamma(A)}_{p\ \text{times}} \ni (\alpha_1, \ldots , \alpha_p) \mapsto
\omega(\alpha_1, \ldots , \alpha_p)\in \Gamma(E), \]
with the differential $\delta_{A, E}:C^{p}(A; E)\rmap C^{p+1}(A; E)$ given by the usual Chevalley-Eilenberg formula:
\begin{eqnarray}
\label{differential2}
\delta_{A, E}(\omega)(\alpha_1, \ldots , \alpha_{p+1}) & = & \sum_{i<j}
(-1)^{i+j}\omega([\alpha_i, \alpha_j], \alpha_1, \ldots , \hat{\alpha_i}, \ldots ,
\hat{\alpha_j}, \ldots \alpha_{p+1}) \nonumber \\
 & + & \sum_{i=1}^{p+1}(-1)^{i+1}
\nabla_{\alpha_i}(\omega(\alpha_1, \ldots, \hat{\alpha_i}, \ldots , \alpha_{p+1})) .
\end{eqnarray}
When $E$ is the trivial line bundle (with $\nabla_{\alpha}= L_{\rho(\alpha)}$, the Lie derivative along $\rho(\alpha)$),
then we omit $E$ from the notation. In particular, the differential on $C^*(A)$ will be denoted $\delta_{A}$.

For instance, if 
$A= TM$ is the tangent bundle of $M$, $C^*(A)= \Omega^*(M)$ and the formula for $\delta_{A}$ becomes
the known Koszul-formula for the De Rham differential. Also, when $A= \mathfrak{g}$ is a Lie algebra (and $M$ consists of one point),  
one recovers the Chevalley-Eilenberg complex computing Lie algebra cohomology with coefficients. 
In general, if $E$ is a representation of the Lie algebroid $A$, then $\Gamma(E)$ becomes a representation of the Lie algebra $\Gamma(A)$, while the Lie algebroid complex  
\[ C^*(A; E)\subset C^*(\Gamma(A); \Gamma(E)) ,\]
is the subcomplex of $C^*(\Gamma(A); \Gamma(E))$ consisting of {\it $C^{\infty}(M)$-multilinear} cochains.

\subsection{Deformation cohomology} It is well-known that deformations of a Lie algebra $\mathfrak{g}$ are controlled by $H^*(\mathfrak{g}, \mathfrak{g})$, the cohomology of $\mathfrak{g}$ 
with coefficients in the adjoint representation, and by the associated differential graded Lie algebra $(C^{*}(\mathfrak{g};\mathfrak{g}), \delta)$ \cite{NR2, NR3}.
Here, the graded Lie algebra structure is the one of Proposition \ref{brackets} applied to the case where $M$ is a point and $E= \mathfrak{g}$, and the differential is the Chevalley-Eilenberg differential, which can also be expressed in terms of the Gerstenhaber bracket with the Lie bracket $m\in C^{2}(\mathfrak{g}, \mathfrak{g})$ of $\mathfrak{g}$ \cite{NR2},  
\[ \delta(c)= [m, c] .\] 
For the case of Lie algebroids, one faces the problem that, in general, Lie algebroids do not have an adjoint representation \cite{coh}, and/or the
$C^{\infty}(M)$-multilinear cochains of the complex $C^*(\Gamma(A); \Gamma(A))$ do not form a sub-complex. However, there is a distinguished subcomplex of $C^{*}(\Gamma(A), \Gamma(A))$ which consists of cocycles which are ``not far from
being $C^{\infty}(M)$-multilinear''. More precisely, since the symbol of an element $D\in Der^*(A)$ is uniquely determined by
$D$, the multiderivations form a subcomplex
\[ Der^{*-1}(A)\subset C^{*}(\Gamma(A); \Gamma(A)) .\]
This is precisely the complex that we have denoted $\de{C}^{*}(A)$ and called {\it the deformation complex
of $A$} at the beginning of this section. That this is 
indeed a subcomplex follows from Proposition \ref{brackets} and the fact that $\delta= [m, -]$, but it can also be shown directly
that, for any multiderivation $D\in Der^{n}(A)$, $\delta(D)$ is again a multiderivation with the symbol
\begin{equation}
\label{lemma-symbol} 
\sigma_{\delta(D)}= \delta(\sigma_{D})+ (-1)^n \rho\circ D,
\end{equation}
where, for $\sigma\in \Gamma(\Lambda^{n}A^{\vee}\otimes TM)$, $\delta(\sigma)\in \Gamma(\Lambda^{n+1}A^{\vee}\otimes TM)$ is given by
\begin{multline}
\delta(\sigma)(\alpha_0,\dots,\alpha_{n})=
\sum_{i} (-1)^{i} [\rho(\alpha_i), \sigma(\alpha_0,\dots,\hat{\alpha_i},\dots,\alpha_{n})]+\\
+\sum_{i<j} (-1)^{i+j}
\sigma([\alpha_i,\alpha_j],\alpha_0,\dots,\hat{\alpha_i},\dots,\hat{\alpha_j},\dots,\alpha_{n}).
\end{multline}
 We summarize the discussion up to this point in the following theorem.

\begin{theorem}
\label{th1-late}
For any algebroid $A$, $\de{C}^{n}(A)= Der^{n-1}(A)$ and $\delta(D)= [m, D]$. In particular, $(\de{C}^{*}(A), \delta)$
is a differential graded Lie algebra (with a shift in degree), and  the deformation cohomology 
$\de{H}^{*}(A)$ is a graded Lie algebra. 
\end{theorem}

\subsection{Alternative descriptions of the deformation complex}
\label{alternative}
There are several different descriptions of the space of multi-derivations $Der^*(E)$ on a vector bundle and of the deformation complex $\de{C}^{*}(A)$ of an algebroid. Here we give a description in terms of derivations (which reveals a connection between our deformation complex and recent unpublished work of D. Roytenberg), while in subsection \ref{lin-cat} we will give a more geometric description in terms of linear multi-vector fields. 
As before, we denote by $E$ a vector bundle over $M$, and we will change the notation to $A$ when dealing with Lie algebroids.

By formulas similar to the classical ones for Lie derivatives of forms along vector fields, any derivation $D: \Gamma(E)\rmap \Gamma(E)$ induces
an $\mathbb{R}$-linear derivation of degree zero on the algebra $C^*(E)= \Gamma(\Lambda^*(E^{\vee}))$ of sections of the exterior bundle (viewed as $C^{\infty}(M)$-valued, $C^{\infty}(M)$-multilinear maps on the powers of $\Gamma(E)$), by
\begin{equation}
\label{Lie-derivative}
L_{D}(c)(s_1, \ldots, s_{q})= L_{\sigma_{D}}(c(s_1, \ldots, s_q))- \sum_{i=1}^{q} c(s_1, \ldots, D(s_i), \ldots , s_q).
\end{equation}
More generally, any $D\in Der^p(E)$ induces a derivation of degree $p$ on $C^*(E)$: if $c\in C^q(E)$, then $L_{D}(c)\in C^{p+q}(E)$ is given by 
\[ L_{D}(c)= (-1)^{pq} \sigma_{D}\circ c- c\circ D,\]
where $c\circ D$ is defined by the Gerstenhaber-type formula (\ref{gert-eq}),  
and, similarly,
\begin{equation} 
\sigma_{D}\circ c (s_1,\dots,s_{p+q})= 
\sum_{\sigma} (-1)^{|\sigma|} L_{\sigma_{D}(s_{\sigma(q+1)} \ldots , s_{\sigma(p+q)})}(c(s_{\sigma(1)}, \ldots , s_{\sigma(q)}))  .\label{nunu}
\end{equation} 
Here $L_{X}$ denotes the Lie derivative along any vector field $X$, and the sum is over all $(q, p)$-shuffles. 
Conversely, any $\mathbb{R}$-linear derivation of degree $p$ on $C^*(E)$ arises in this way, since it is uniquely determined by what it does on $C^{\infty}(M)$ and $\Gamma(E)$. For later reference we give the explicit formulas for $L_{D}$ applied to $f\in C^0(E)= C^{\infty}(M)$ and to  $\xi\in C^1(E)= \Gamma(E^{\vee})$ , when $D\in Der^{p-1}(E)$:
\begin{equation}
\label{low0}
L_{D}(f)(s_1, \ldots , s_{p-1}) = \sigma_{D}(s_1, \ldots , s_{p-1}) (f),
\end{equation}
\begin{equation}
\label{low1}
L_{D}(\xi)(s_1, \ldots , s_{p})=\sum_{i=1}^{p} (-1)^{p-i} L_{\sigma_D(s_1, \ldots \widehat{s_i} \ldots , s_p)}(\xi(s_i))- \xi(D(s_1, \ldots , s_p)).
\end{equation}
The conclusion is that $Der^{*}(E)$ is isomorphic to the algebra of derivations of $C^{*}(E)$ (as graded Lie algebras!).  With this, Lemma \ref{rer-alg} translates into
the following well-known observation \cite{KM}.

\begin{corollary} Given a vector bundle $E$, there is a 1-1 correspondence between Lie algebroid structures on $E$ and derivations $\delta$ of degree $1$ on the algebra $C^*(E)= \Gamma(\Lambda E^{\vee})$, satisfying $\delta^2= 0$. 
\end{corollary}

Note that, if $E= A$ is a Lie algebroid, we obtain an action of $\de{C}^{*}(A)$ on $C^{*}(A)$ (be aware of the degree shift!), and the differential $\delta_{A}$ of $C^{*}(A)$ coincides with $L_{m}$, where $m\in \de{C}^{2}(A)$ is the Lie bracket of $A$. The discussion above and a careful computation of the boundaries shows that, conversely, one recovers the differential graded Lie algebra $(\de{C}^{*}(A), \delta)$ as the algebra of derivations of the differential graded algebra $(C^*(A), \delta_A)$. We will come back to this point is subsection \ref{act-Lie-coh} below.

\section{Deformation cohomology and deformations}%
\label{Deformations}              %

In this section we examine our deformation cohomology of Lie algebroids in low degrees. In degree zero,
\[ \de{H}^{0}(A)= Z(\Gamma(A)), \]
the center of the infinite-dimensional Lie algebra of sections of $A$. In degrees 1 and 2 we will show that
the usual interpretations for Lie algebras extend to the context of Lie algebroids. In particular, degree 2 
cohomology classes will be seen to correspond to deformations, thus justifying the name ``deformation cohomology''.

\subsection{$\de{H}^{1}$ and derivations}
\label{sec:1st cohom}
Recall that a {\it derivation of a Lie algebroid $A$} is a linear map
$D: \Gamma(A)\rmap \Gamma(A)$  which is both a vector bundle derivation (see subsection \ref{Multiderivations}) and a derivation
with respect to the Lie bracket:
\[ D([\alpha, \beta])= [D(\alpha), \beta]+ [\alpha, D(\beta)] .\]
These derivations form a Lie algebra $Der(A)$ under the commutator bracket. This Lie algebra was studied
in \cite{MX, MM}. The {\it inner derivations}, i.e. those 
of the form $[\alpha, -]$, form an ideal in $Der(A)$, and the quotient, denoted $OutDer(A)$, is the Lie algebra of
{\it outer derivations} of $A$. It is immediate from the definitions that 
\[ \de{H}^{1}(A)= OutDer(A).\]

This space can also be interpreted as the Lie algebra of the (infinite-dimensional) group of outer automorphisms of $A$. 
For later use, we make this statement more precise. 
The passage from the
infinitesimal side (derivations) to the global side (automorphisms) is
via flows of derivations (see the Appendix in \cite{CrFe1}). Given
$D\in Der(E)$, the flow $\Phi_{D}^{t}$ of $D$ is a 1-parameter group
of bundle isomorphisms of $E$, covering the flow $\phi_{\sigma_D}^{t}$ of
the symbol $\sigma_D$ of $D$:
\[ \Phi_{D}^{t}(x): E_{x}\to E_{\phi_{\sigma_D}^{t}(x)}.\]
It is characterized uniquely by the property
\begin{equation}
\label{eq:flow}
{\left. \frac{d}{dt}\right|}_{t= 0} (\Phi_{D}^{t})^{*}\beta =D(\beta),
\qquad \forall \beta\in \Gamma(E).
\end{equation}
In this relation, $(\Phi_{D}^{t})^{*}(\beta)(x) =
\Phi_{D}^{t}(\beta(\phi_{\sigma_D}^{-t}(x))$. In general, $\Phi_{D}^{t}(x)$
will be defined whenever $\phi_{\sigma_D}^{t}(x)$ is, and one has to deal
with local bundle maps, defined only over some open sets. Since this is
standard, and not relevant for the present discussion, we will assume
that the vector fields $\sigma_D$ are complete.

\begin{lemma} 
Let $A$ be a Lie algebroid. A derivation $D\in Der(A)$ is a Lie
algebroid derivation if and only if the maps $\Phi_{D}^{t}$ are Lie algebroid
automorphisms.
\end{lemma}

The proof is standard. In the case of an inner derivation
$D=\ad_{\alpha}$, one denotes $\Phi_{D}^{t}$ by $\Phi_{\alpha}^{t}$, and
calls it the (infinitesimal) flow of $\alpha$. These flows play an
essential role in \cite{CrFe1}, and they can be thought of as the
inner automorphisms of $A$. Therefore, we can think of the
Lie algebra $\de{H}^{1}(A)$ as the Lie algebra of the (infinite
dimensional) group of outer automorphisms of $A$.

\subsection{$\de{H}^{2}$ and deformations}
\label{sec:2nd cohom}
We now explain the relevance of the deformation complex to
deformations of Lie algebroids.

\begin{definition}
Let $A$ be a fixed vector bundle, and $I\subset \mathbb{R}$ an interval. 
\begin{enumerate}
\item[(i)] A {\it family} of Lie algebroids over $I$ is a
  collection $(A_t)_{t\in I}$ of Lie algebroids
  $A_t=(A,[\cdot,\cdot]_t,\rho_t)$ varying smoothly with respect to $t$;
\item[(ii)] We say that $(A_t)_{t\in I}$ and $(A'_t)_{t\in I}$ are
  {\it equivalent families} of Lie algebroids if there exists a
  family of Lie algebroid isomorphisms $h_{t}: A_{t}\to A'_{t}$,
  depending smoothly on $t$;
\item[(iii)] A {\it deformation} of a Lie algebroid
  $(A,[\cdot,\cdot],\rho)$ is a family $(A_t)_{t\in I}$ of Lie
  algebroids over an interval containing the origin with $A_0=A$;
\item[(iv)] Two deformations $(A_t)_{t\in I}$ and $(A'_t)_{t\in I}$ of
  a Lie algebroid $A$ are said to be {\it equivalent deformations}
  if there exists an equivalence $h_{t}: A_{t}\to A'_{t}$ with
  $h(0)=Id$.
\end{enumerate} 
\end{definition}

We will say that a family (or deformation) is {\it trivial} if it is
equivalent to the constant family (deformation). We have the following interpretation for the elements of the second
cohomology group $\de{H}^2(A)$ in terms of deformations.

\begin{proposition}
\label{prop:deformation}  
Let $A_{t}= (A, [\cdot, \cdot]_{t}, \rho_{t})$ be a deformation of the
Lie algebroid $A$. Then
\[ c_{0}(\alpha,\beta)= \left.\frac{d}{dt} [\alpha, \beta]_{t} \right|_{t=0}\]
defines a cocycle $c_{0}\in \de{C}^{2}(A)$. The corresponding
cohomology class in $\de{H}^{2}(A)$ only depends on the equivalence
class of the deformation. 
\end{proposition}

\begin{proof}
Let us denote, as before, by $m_{t}\in Der^1(A)$ the Lie bracket $[\cdot,
\cdot]_{t}$. Since $Der^1(A)= \Gamma(\mathcal{D}^1(A))$, $c_0$ is a multiderivation
(its symbol is $\sigma= \frac{d}{d t}\rho_t|_{t=0}$). 
Taking derivatives at $t=0$ in the equation 
$[m_t, m_t]=0$ we obtain $[c_{0},m]= 0$, i.e., $\delta(c_0)=0$ so that $c_0$ is a cocycle.

Assume now that $A'_{t}$ is another deformation, and denote the
associated class by $c'_{0}$. Assume also that $h_{t}$ defines an
equivalence between $A_{t}$ and $A'_{t}$. We use the same notation
$h_{t}$ for the map induced at the level of sections, and we consider
the derivation $D$ defined by:
\[ D= {\left.\frac{d}{dt}\right|}_{t=0} h_{t}: \Gamma(A)\to \Gamma(A).\]
Since $h_{t}$ is a Lie algebroid map we have
\[ h_{t}([\alpha,\beta]_{t})= [h_t(\alpha), h_t(\beta)]'_{t} .\]
Taking derivatives of both sides and setting $t=0$, we obtain
\[ 
D([\alpha,\beta])+ c_{0}(\alpha, \beta)= c'_{0}(\alpha,\beta)+[D(\alpha),\beta]+[\alpha,D(\beta)],
\]
which means precisely that $c_{0}-c'_{0}=\delta(D)$.
\end{proof}

\begin{remark}\label{var-bundle} Note that, when talking about deformations, one can also allow the vector bundle $A$ itself to vary smoothly with respect to $t$ 
(in the sense that the $A_t$'s together fit into a smooth vector bundle over $M\times \mathbb{R}$). Indeed, to any such family $A_t$ of algebroids, one can associate (an equivalence class of) a family of algebroids with constant vector bundle $A= A_0$, by choosing vector bundle isomorphisms $\phi_t: A_t\rmap A$ (which is possible because the real line is contractible) and transporting the bracket of $A_t$ to a bracket $[\cdot, \cdot]_t$ on $A$. 
Of course, different choices of $\phi_t$ produce equivalent deformations in the sense above. In particular, the cohomology class of the deformation is defined unambiguously for any deformation of $A$ with possible varying vector bundle.
\end{remark} 

If $(A_t)_{t\in I}$ is a trivial deformation of a Lie algebroid $A$
then obviously we must have $[c_0]=0$. The converse is not
true, but our next result gives a partial converse.

\begin{theorem} 
\label{rigidity}
Let $(A_{t})_{t\in I}=(A,[\cdot,\cdot]_{t},\rho_t)$ be a family of Lie
algebroids. Then
\[ c_t(\alpha,\beta)= \frac{d}{dt} [\alpha, \beta]_{t}.\]
defines a cocycle $c_{t}\in \de{C}^{2}(A_t)$. 
If $M$ is compact, then the following are equivalent
\begin{enumerate}
\item[(i)] The family $(A_{t})_{t\in I}$ is
trivial.
\item[(ii)] The classes $[c_{t}]\in \de{H}^{2}(A_t)$ vanish
smoothly with respect to $t$, i.e., $c_{t}= \delta(D_t)$ for a smooth
family $D_{t}$ of $1$-cochains. 
\end{enumerate}
\end{theorem}

\begin{proof} 
The fact that, for each $t$, $c_{t}$ is a cocycle and its cohomology
class depends only on the equivalence class of $(A_{t})_{t\in I}$
follows from Proposition \ref{prop:deformation}, since we can view
$A_s$ as a deformation of $A_t$.

Let $(A_{t})_{t\in I}$ be a trivial family, so that there exist Lie
algebroid isomorphisms $h_t:A_t\to A$, to a fixed Lie algebroid
$A$. We define $D_t:\Gamma(A_t)\to \Gamma(A_t)$ by letting
\[ h_t(D_t(\alpha))=\frac{d}{dt} h_{t}(\alpha).\]
Differentiating both sides of
\[ h_t([\alpha,\beta]_t)=[h_t(\alpha),h_t(\beta)],\]
we obtain:
\begin{align*} 
h_t(D_t([\alpha,\beta]_t)+c_t(\alpha,\beta))&=[h_t(D_t(\alpha)),h_t(\beta)]+
[h_t(\alpha),h_t(D_t(\beta))]\\
&=h_t([D_t(\alpha),\beta]_t+[\alpha,D_t(\beta)]_t).
\end{align*}
This shows that
\[ c_t(\alpha,\beta)=[D_t(\alpha),\beta]_t+[\alpha,D_t(\beta)]_t-D_t([\alpha,\beta]_t),\]
which means precisely that $c_t=\delta(D_t)$, so the classes
$[c_{t}]\in \de{H}^{2}(A_t)$ vanish smoothly with respect to $t$.

Conversely, suppose that $(A_{t})_{t\in I}$ is a family of Lie
algebroids, such that $c_{t}= \delta (D_t)$ for a smooth family $D_{t}$ of
$1$-cochains. Each $D_t$ is a derivation, and since $M$ is compact, the flow
$\Phi^s_{D_t}$ of $D_t$ is defined for all $s$, for each fixed
$t$. Denote by $h_t$ the flow at time $t$. From the defining relation
(\ref{eq:flow}) for the flow of a derivation, we have:
\[ h_t(D_t(\alpha))=\frac{d}{dt} h_{t}(\alpha).\]
We claim that
\begin{equation}
\label{eq:equivalence}
h_t([\alpha,\beta]_t)=[h_t(\alpha),h_t(\beta)]_0.
\end{equation}
This shows that $h_t:A_t\to A_0$ gives an equivalence to a constant
family, so the family $(A_{t})_{t\in I}$ will be trivial.

To prove the claim, we just observe that (\ref{eq:equivalence}) holds
at $t=0$, and that the derivative of both sides of
(\ref{eq:equivalence}) are equal. In fact, as we saw above, the
derivative is precisely the condition $c_{t}= \delta(D_t)$. Therefore
equality holds for all $t$.
\end{proof}

\section{Relations to known cohomologies and particular cases}         %
\label{Particular cases}          %

In this section we look at some particular classes of algebroids, and we relate $\de{H}^{*}(A)$ to 
known cohomology theories. We begin by simply mentioning the following two extreme cases.

\subsection{Lie algebras} Since the Lie algebra case was partially used as inspiration for our constructions, it is clear that the deformation complex $\de{C}^{*}(\mathfrak{g})$ of a Lie algebra $\mathfrak{g}$ is the usual Chevalley-Eilenberg complex $C^{*}(\mathfrak{g}, \mathfrak{g})$ with coefficients in the adjoint representation, and one recovers the classical relation between Lie algebra deformations, the graded Lie algebra $C^{*}(\mathfrak{g}, \mathfrak{g})$, and the cohomology groups $H^{*}(\mathfrak{g}, \mathfrak{g})$ \cite{NR2, NR3}.

\subsection{Tangent bundles} 
Let us now look at the case where $A= TM$ is the tangent bundle of a manifold $M$. The first remark is that closed cocycles correspond to vector valued forms on $M$:
\[ Z^k\de{C}^{*}(TM)\cong \Gamma(\Lambda^{k-1}T^{\vee}M\otimes TM) ,\]
where $T^{\vee}M$ denotes the cotangent bundle. 
This follows from (\ref{lemma-symbol}) which shows that any cocycle $D$ in the deformation complex is determined by its symbol $\sigma_{D}$ by 
$D= (-1)^k\delta(\sigma_D)$. Actually, the same formula shows that the map which associates to $D$ its symbol,
viewed as a map of degree $-1$ on $\de{C}^{*}(TM)$, defines a homotopy in the deformation complex (up to a sign). Hence:

\begin{corollary} 
Any closed multiderivation on $TM$ is exact (i.e. $\de{H}^{*}(TM)= 0$).
\end{corollary}

A different (and more general) argument will be presented in subsection \ref{reg-case} below.

\subsection{Cotangent bundles} We now consider the graded Lie algebra of multiderivations on the cotangent bundle $T^{\vee}M$ of a manifold $M$,
and we relate it to the known graded Lie algebra of multi-vector fields on $TM$, $(\mathcal{X}^{*}(M), [\cdot, \cdot])$.
Recall that $\mathcal{X}^{p}(M)=
\Lambda^p(T^{\vee}M)$, where the Lie algebra degree of a $p$-vector field is $(p-1)$, so that the graded antisymmetry reads:
\[ [X, Y]= - (-1)^{(p-1)(q-1)}[Y, X] \]
for $X\in \mathcal{X}^{p}(M)$, $Y\in \mathcal{X}^{q}(M)$. Recall also that $[\cdot, \cdot]$ is the Nijenhuis-Schouten bracket; explicitly, 
for $X_i, Y_j\in \mathcal{X}(M)$ and $f\in C^{\infty}(M)$,
\[ [X_1\wedge \ldots \wedge X_p, Y_1\wedge \ldots \wedge Y_q]= \sum_{i, j} (-1)^{i+j} [X_i, Y_j]\wedge X_1\wedge \ldots \widehat{X_i}\ldots \wedge X_p \wedge Y_1\wedge \ldots \widehat{Y_j}\ldots \wedge Y_q ,\]
\[ [X_1\wedge \ldots \wedge X_p, f]= \sum_{i=1}^{n} (-1)^{p-i} L_{X_i}(f) X_1\wedge \ldots \widehat{X_i} \ldots \wedge X_p .\]

To state the next proposition, we also need the following notation:  
for $X\in \mathcal{X}^n(M)$, we view $X$ as an antisymmetric map depending on $n$ one-forms $\omega_1, \ldots , \omega_n$. Fixing the 
first $n-1$ of them defines a linear map $X(\omega_1, \ldots , \omega_{n-1}, -)$ on $T^{\vee}M$, hence a vector field. We denote this vector field 
by $X^{\sharp}(\omega_1, \ldots , \omega_{n-1})$, so that $X^{\sharp}$ becomes a linear antisymmetric map 
\[ X^{\sharp} : \underbrace{T^{\vee}M\otimes \ldots \otimes T^{\vee}M}_{(n-1)\ \text{times}}\rmap TM .\]

\begin{proposition}
\label{prop-ctg} For any $X\in \mathcal{X}^n(M)$, there exists an unique $D_X\in Der^{n-1}(T^{\vee}M)$ with symbol
$X^{\sharp}$ and satisfying 
\[ D_X(df_1, \ldots , df_n)= d( X(df_1, \ldots , df_n)),\]
for all $f_i\in C^{\infty}(M)$. Explicitly, for all $\omega_i\in \Omega^1(M)$, $1\leq i\leq n$, 
\[ D_X(\omega_1, \ldots , \omega_n)= \sum_{i=1}^{n} (-1)^{n-i} L_{X^{\sharp}(\omega_1, \ldots , \widehat{\omega_i}, \ldots , \omega_n)}(\omega_i)- (n-1) d(X(\omega_1, \ldots , \omega_n)) ,\]
Moreover, the map $\mathcal{X}(M)\rmap \Der(T^{\vee}M)$, $X\mapsto D_X$
is a map of graded Lie algebras.
\end{proposition}

\begin{proof} The uniqueness part is clear. That the explicit formula for $D_X$ defines a multi-derivation with the desired 
properties follows by direct computation. 
To prove the last part, we 
recall first that $(\mathcal{X}^{*}(M), \wedge, [\cdot, \cdot])$ is a Gerstenhaber algebra, where $\wedge$ is the exterior product. 
Apart from the fact that $(\mathcal{X}^{*}(M), [\cdot, \cdot])$ is a graded Lie algebra (with the Lie algebra degree of $X\in \mathcal{X}^p(M)$ equal to $p-1$), this also means that $(\mathcal{X}^{*}(M), \wedge)$ is a graded algebra (without degree-shift!), while the two structures are related by the Leibniz rule:
\[ [X, Y\wedge Z]= [X, Y]\wedge Z+ (-1)^{(p-1)q} Y\wedge [X, Z],\]
for all $X\in \mathcal{X}^{p}(M)$, $Y\in \mathcal{X}^{q}(M)$, $Z\in \mathcal{X}^{r}(M)$. More formally, this equation means that 
the map $X\mapsto [X, -]$ is a graded map from $(\mathcal{X}^*(M), [\cdot, \cdot])$ into the graded Lie algebra 
$Der(\mathcal{X}^*(M), \wedge)$ of derivations on the graded algebra $(\mathcal{X}^*(M), \wedge)$.

On the other hand, by \ref{alternative} applied to $T^{\vee}M$, we have an isomorphism $D\mapsto L_{D}$ from the graded Lie algebra
$Der(T^{\vee}M)$ of multiderivations on the cotangent bundle into the graded Lie algebra of derivations on $(\mathcal{X}^*(M), \wedge)$.
We claim that $L_{D_X}= [X, -]$. Since this is an equality of derivations on $(\mathcal{X}^*(M), \wedge)$, it suffices to show that
$L_{D_X}$ and $[X, -]$ are equal on functions and on vector fields. This follows again by direct computation, using the formulas for the Nijenhuis-Schouten bracket given above, and the formulas (\ref{low0}), (\ref{low1}) describing $L_{D}$ in low degrees for $D= D_X$.
\end{proof}

\subsection{Poisson manifolds I} Here we describe the relation of our deformation cohomology with Poisson cohomology.
Recall that a Poisson manifold is a pair $(P, \pi)$ where $\pi \in \Gamma(\wedge^2(TP))$
is a bivector with the property that the ``Poisson bracket'' $\{f, g\}= \pi(df, dg)$
on $C^{\infty}(P)$ satisfies the Jacobi identity, or, equivalently, the Nijenhuis-Schouten bracket $[\pi, \pi]$ vanishes.
It is well-known that a Poisson structure on $P$ induces an algebroid structure on $T^{\vee}P$ with anchor $\pi^{\sharp}$ (defined by $\beta(\pi^{\sharp}(\alpha))= \pi(\alpha, \beta)$, for all 1-forms $\alpha$ and $\beta$).  The bracket is 
usually introduced either by the explicit formula
\[ [\alpha, \beta]= L_{\pi^{\sharp}(\alpha)}(\beta)- L_{\pi^{\sharp}(\beta)}(\alpha)- d(\pi(\alpha, \beta)),\]
or by saying that it is the unique Lie algebroid structure on $T^{\vee}P$ with anchor $\pi^{\sharp}$ and the property that $[df, dg]= d\{f, g\}$.

It is interesting to relate this to the previous proposition. Since $[\pi, \pi]= 0$, the induced derivation $D_{\pi}\in Der^1(T^{\vee}P)$ satisfies 
the same formula, hence it defines an algebroid structure on $T^{\vee}P$ by Lemma \ref{rer-alg}. This coincides with the known algebroid structure, and the two ways of describing the structure correspond to the descriptions of $D_{\pi}$ in Proposition
\ref{prop-ctg}.

Next, given the Poisson manifold $P$, the De Rham cohomology of the induced Lie algebroid $T^{\vee}P$ (with coefficients in the trivial line bundle) is  known as the Poisson cohomology of
$P$, denoted $H^{*}_{\pi}(P)$. The defining complex is the complex of multivector fields
\[ C^{*}_{\pi}(P)= \Gamma(\wedge^{*}TP), \]
with boundary $d(X)= [\pi, X]$, where $[\cdot, \cdot]$ is the Nijenhuis-Schouten bracket (see e.g. \cite{Va} and the references therein). In particular, $C^{*}_{\pi}(P)$ is a
dgla.

Finally, if $\pi_{t}$ is a family of Poisson structures on $P$ with $\pi_0= \pi$, then 
taking derivatives with respect to $t$ in $[\pi_t, \pi_t]= 0$ at $t=0$, we see that 
\begin{equation}
\label{coh-P} 
[\frac{d}{dt}|_{t=0} \pi_{t}]\in H^{2}_{\pi}(P) 
\end{equation}
is a well-defined cohomology class. This is known as the cohomology class associated to the deformation
$\pi_{t}$. 
Using the last part of Proposition \ref{prop-ctg} and the remark above that $D_{\pi}$ is the Lie bracket of $T^{\vee}P$, we deduce:

\begin{corollary} For any Poisson manifold $(P, \pi)$, the map $X\mapsto D_X$ makes the Poisson complex $C^{*}_{\pi}(P)$ into 
a dg Lie sub-algebra
of the deformation complex $\de{C}^{*}(T^{\vee}P)$. In particular, there is an induced map of graded Lie algebras
\[ i: H^{*}_{\pi}(P)\rmap \de{H}^{*}(T^{\vee}P) .\]
Moreover, if $\{\pi_t\}$ is a deformation of $\pi$, then $[\cdot, \cdot]_{\pi_t}$ defines a deformation
of the Lie algebroid $(T^{\vee}P, [\cdot, \cdot]_{\pi})$, and the associated cohomology classes (i.e. (\ref{coh-P}),
and the one of Proposition \ref{prop:deformation}, respectively) are related by the map $i$.
\end{corollary}

\subsection{Foliations} 
We now look at the particular case of regular foliations on a manifold $M$, i.e. sub-bundles $\mathcal{F}\subset TM$
(of vectors tangent to the leaves) which are involutive (i.e. $[\Gamma(\mathcal{F}), \Gamma(\mathcal{F})]\subset \Gamma(\mathcal{F})$).
Regular foliations are the same thing as Lie algebroids with injective anchor map. Similar to the adjoint representation
of a Lie algebra, any foliation $\mathcal{F}$ has a canonical representation on the normal bundle $\nu= TM/\mathcal{F}$, called the
Bott representation, defined by $\nabla_{X}(\overline{Y})= \overline{[X, Y]}$. The resulting cohomology $H^*(\mathcal{F}; \nu)$ is known as the foliated
(or leafwise) cohomology with coefficients in the normal bundle, and was investigated by Heitsch \cite{Hei} in connection 
with deformations of foliations. Explicitly, given a family of foliations $\mathcal{F}_t$ with $\mathcal{F}_0= \mathcal{F}$, Heitsch defines 
\begin{equation}
\label{H-class} 
c_{0}(v)= \pi_{0}^{\perp}(\frac{d}{dt}|_{t=0} \pi_{t}(v)) ,
\end{equation}
for all $v\in \Gamma(\mathcal{F})$, where $\pi_{t}:TM\rmap \mathcal{F}_{t}$ and $\pi_{t}^{\perp}:TM\rmap \nu_{t}$ are the orthogonal projections
with respect to a Riemannian metric. Then $c_{0}$ defines a cohomology class 
\[ [c_{0}]\in H^{1}(\mathcal{F}; \nu), \]
independent of the metric, which we will call {\it the Heitsch characteristic class} of the deformation.

On the other hand, each such deformation defines a deformation of $\mathcal{F}$ viewed as an algebroid, hence a class in $\de{H}^{2}(\mathcal{F})$. We have:

\begin{proposition} For any foliation  $\mathcal{F}$,
\[ \de{H}^{*}(\mathcal{F})\cong H^{*-1}(\mathcal{F}; \nu) \ .\]
Moreover, for a deformation $\mathcal{F}_t$ of $\mathcal{F}$, the induced cohomology 
class in $\de{H}^{2}(\mathcal{F})$ corresponds to the Heitsch characteristic class of the deformation.
\end{proposition}

Note that, in degree zero,
\[ l(M, \mathcal{F}):= H^{0}(\mathcal{F}; \nu) \]
is well known in foliation theory as the Lie algebra of transversely projectable vector fields. 
It consists of sections $\overline{X}$ of the normal bundle with the property that
$[X, \Gamma(\mathcal{F})\subset \Gamma(\mathcal{F})$.  The previous proposition together with Theorem \ref{th1-late}
implies that
the Lie algebra structure on $l(M, \mathcal{F})$ is just the degree zero part of a graded Lie algebra
structure on $H^{*}(\mathcal{F}; \nu)$.

\begin{proof} The first part of the proposition will be immediate from Theorem \ref{regular} below. Here we sketch a different argument,
similar to the one in the case of tangent bundles. Note first that equation (\ref{lemma-symbol}) implies that if $D\in \de{C}^{k}(\mathcal{F})$ is a cocycle, then
$D$ is uniquely determined by its symbol $\sigma_{D}\in C^{k-1}(\mathcal{F}; TM)$. On the other hand, the symbol of $D$ projects to an element $\overline{\sigma}_{D}\in C^{k-1}(\mathcal{F}; \nu)$
which will be a cocycle. One then checks directly that the induced cohomology class only depends on the cohomology class of $D$, 
and this defines the desired isomorphism. (Alternatively, one can easily chase the sequences appearing in the proof of Theorem \ref{regular}).

Let us now go to the second part. We are in the case described in Remark \ref{var-bundle} where the vector bundle is varying also,
hence we first need to trivialize $\mathcal{F}_t$ as a family of vector bundles, and then consider the induced brackets $[\cdot, \cdot]_t$ on $\mathcal{F}$. The trivialization we will use is the one induced by the parallel transport (with respect to $t$) of the canonical partial connection $\nabla_{\frac{d}{dt}}= \pi_{t}\circ \frac{d}{dt}$ on $(\mathcal{F}_{t})_t$ viewed as a bundle over $M\times \mathbb{R}$. Since the induced trivialization $T_t: \mathcal{F}\rmap \mathcal{F}_t$ sends $v$ into its parallel transport $T_t(v)$ at time $t$, it is the solution of the equation:
\[ \pi_{t}(\frac{d}{dt} T_t(v)) = 0, T_{0}(v)= v .\]
The deformation cocycle induced by the deformation is $c= \frac{d}{dt}|_{t=0}[\cdot , \cdot]_t$. From the first part,
we know that the cocycle $\sigma\in C^{1}(\mathcal{F}; \nu)$ corresponding to $c$ is the symbol of $c$ (i.e. the anchor $T_t$ of
$[\cdot , \cdot]_t$) projected down to $\nu$:
\[ \sigma(v)= \pi_{0}^{\perp}(\frac{d}{dt}|_{t= 0} T_t(v)) \in \Gamma(\nu).\]
We claim that this coincides with $c_0$ of (\ref{H-class}) defining the Heitsch class. Since both $c_0(v)$ and $\sigma(v)$ are tangent to $\mathcal{F}$, we only have to show that the orthogonal projection $\pi_{0}^{\perp}$ onto $\mathcal{F}$ kills their difference.
This follows immediately by taking derivatives in $\pi_{t}^{\perp}(T_t(v)- \pi_t(v))= 0$ at $t= 0$. 
\end{proof}

\subsection{Lie algebra actions on manifolds}
Given an action of a Lie algebra $\mathfrak{g}$ on a manifold $M$, i.e. a Lie algebra map $\rho_M: \mathfrak{g}\rmap \mathcal{X}(M)$, one has an induced Lie algebroid structure on $\mathfrak{g}\times M$ (the trivial vector bundle over $M$ with fiber $\mathfrak{g}$): the anchor is the map $\rho$ defining the action, while the bracket is determined by its values on constant sections (the bracket between the constant sections $v, w\in \mathfrak{g}\subset \Gamma(\mathfrak{g}\times M)$ is the constant section $[v, w]$), and the Leibniz rule. The resulting algebroid is denoted by $\mathfrak{g}\ltimes M$ and is called {\it the action algebroid} associated to the infinitesimal action of $\mathfrak{g}$ on $M$. Such algebroids (and the study of their deformations) are important in the linearization problem of Poisson manifolds, although the case where $\mathfrak{g}$ is 
1-dimensional is already interesting (the action on $M$ will be described by a vector field $X$ on $M$, and the integration of the associated action algebroid is nothing but the integration of the vector field $X$).  

Such an action algebroid $A= \mathfrak{g}\ltimes M$ has two canonical representations. First of all, there is a natural action on the trivial vector bundle over $M$ with fiber $\mathfrak{g}$, $\mathfrak{g}_{M}= \mathfrak{g}\times M$, given by the unique $A$-connection $\nabla$ with the property that $\nabla_v(w)= [v, w]$, the bracket of $\mathfrak{g}$, for all $v\in \mathfrak{g}$ (viewed as constant sections of $A$) and $w\in \mathfrak{g}$ (viewed as constant sections of the representation). 
Note that the De Rham cohomology of $A$ with coefficients in this representation is
\[ H^*(A; \mathfrak{g}_{M})= H^{*}(\mathfrak{g}; \Gamma(\mathfrak{g}_{M})),\]
the cohomology of the Lie algebra $\mathfrak{g}$ with coefficients in $\Gamma(\mathfrak{g}_{M})= C^{\infty}(M; \mathfrak{g})$.

Next, there is a similar action of $A$ on $TM$, where the connection is determined by $\nabla_v(X)= [\rho(v), X]$, for $v\in \mathfrak{g}$. The resulting cohomology is isomorphic to the cohomology of $\mathfrak{g}$ with coefficients in $\mathcal{X}(M)$:
\[ H^*(A; TM)= H^{*}(\mathfrak{g}; \mathcal{X}(M)).\]

\begin{proposition}
\label{lie-acti} For any action of a Lie algebra $\mathfrak{g}$ on a manifold $M$, the deformation cohomology of the associated action algebroid
 $A= \mathfrak{g}\ltimes M$ fits into a long exact sequence
\[ \ldots \rmap H^{n-1}(A; TM)\rmap \de{H}^{n}(A)\rmap H^{n}(A; \mathfrak{g}_{M}) \stackrel{\delta}{\rmap} H^{n}(A; TM)\rmap \ldots .\]
\end{proposition}

\begin{proof} We have a short exact sequence 
\[0\rmap  C^{n-1}(A; TM)\rmap \de{C}^{n}(A)\stackrel{\pi}{\rmap} C^{n}(A; \mathfrak{g}_{M})\rmap 0 ,\]
where $\pi$ associates to $D\in \de{C}^{n}(A)$ the unique $\pi(D)= c\in C^n(A; \mathfrak{g}_{M})$ such that
$D(v_1, \ldots , v_n)= c(v_1, \ldots , v_n)$ on constant sections.
\end{proof}

\subsection{The regular case} 
\label{reg-case}
We now relate the deformation cohomology of $A$ to the De Rham cohomology of $A$ with coefficients, in the case
where $A$ is regular, i.e. when $\rho$ has constant rank. In this case the image of $\rho$ defines a (regular) foliation $\mathcal{F}$ of $M$, and the isotropy bundle and the normal bundle:
\[ \mathfrak{g}(A)= Ker(\rho), \ \ \nu= TM/\mathcal{F} \]
are both representations of $A$: the first with $\nabla_{\alpha}(\beta)= [\alpha, \beta]$, and the second one with the Bott 
connection $\nabla_{\alpha}(\overline{X})= \overline{[\rho(\alpha), X]}$. With these notations, we will show that:

\begin{theorem} 
\label{regular}
For any regular Lie algebroid $A$,
there is an associated long exact sequence
\[ \ldots \rmap H^n(A; \mathfrak{g}(A))\rmap \de{H}^{n}(A)\rmap H^{n-1}(A; \nu) \stackrel{\delta}{\rmap} H^{n+1}(A; \mathfrak{g}(A))\rmap \ldots ,\]
where, as above, $\mathfrak{g}$ is the isotropy Lie algebra of $A$, and $\nu$ is the normal bundle of the foliation induced by $A$.
\end{theorem}

\begin{proof} To prove the theorem, we will introduce two auxiliary complexes $C_{1}^{*}$ and $C_{2}^{*}$ 
which fit into exact sequences of cochain complexes:
\begin{eqnarray}
0\rmap C^{*}(A; \mathfrak{g}(A))\stackrel{j}{\rmap} \de{C}^{*}(A)\stackrel{\rho}{\rmap} C_{1}^{*}\rmap 0 ,\label{seq1} \\
0\rmap  C_{1}^{*} \stackrel{i}{\rmap}  C_{2}^{*} \stackrel{\pi}{\rmap} C^{*}(A; \nu)\rmap 0 \label{seq2}
\end{eqnarray}
and such that $C_{2}^{*}$ has zero cohomology. Then the long exact sequence associated to (\ref{seq2}) 
implies that $H^{n}(C_{1}^{*})\cong H^{n-1}(A; \nu)$, which plugged into the long exact sequence associated 
to (\ref{seq1}) will prove our theorem.

First of all, $C_{2}^{n}$ consists of antisymmetric
multilinear maps
\[ D: \underbrace{\Gamma(A) \otimes \ldots \otimes\Gamma(A)}_{n\ \text{times}}\rmap
\Gamma(TM) \] 
together with $\sigma_{D}\in \Gamma(\wedge^{n-1}A^{\vee}\otimes TM)$ (the symbol of $D$), such that
\[ D(fv_1, v_2, \ldots , fv_n)= f D(v_1, v_2, \ldots , v_n)+ \sigma_{D}(v_1, \ldots ,
v_{n-1})(f) \rho(v_n) \] 
for all functions $f$, and all sections $v_i$ (when $\rho\neq 0$ then, like for multiderivations, $\sigma_D$ is uniquely determined by $D$).
Note the similarity with $\de{C}^{*}(A)$, and we complete this similarity by defining the differential 
$\delta$ on $C^{*}_{2}$ by the same formula as in (\ref{differential}) and by equation (\ref{lemma-symbol}).
Next, $C_{1}^{*}$ is defined as the subcomplex of $C_{2}^{*}$ consisting of those $D$ which are $\Gamma(\mathcal{F})$-valued, and the maps in the exact sequences above are the obvious ones. The exactness of the sequences is evident, 
except maybe for the right hand side of (\ref{seq1}) that we now explain. Let $D\in C_{1}^{*}$. Due to the exactness
of the sequence in Lemma \ref{deriv}, we find $D^{'}\in \de{C}^{n}(A)$ whose symbol $\sigma_{D'}$ coincides with $\sigma_{D}$. 
Then $D- \rho(D^{'})$ is multilinear, i.e. comes from a vector bundle map
$\wedge^nA\rmap \mathcal{F}$. Since $\rho: A\rmap \mathcal{F}$ is surjective, we find a map
$a:\wedge^nA\rmap TM$ such that $D- \rho(D^{'})= \rho(a)$, hence $D= \rho(D{'}+ a)$. 

To prove that $C_{2}^{*}$ is acyclic, we remark that any $a\in \Gamma(\wedge^{n-1}A^{\vee}\otimes TM)$
can be viewed as an element in $C^{n-1}_{2}$ with zero symbol, and $\delta(a)\in C^{n}_{2}$ has as symbol
$\sigma_{\delta(a)}= (-1)^n a$. Assume now that $D\in C^{n}_{2}$ is a cocycle. Then 
$D^{'}= D+ (-1)^{n-1}\delta(\sigma_{D})\in C^{n}_{2}$
will have $\sigma_{D'}= 0$, by the previous formula applied to $a= \sigma_{D}$. Hence we can apply the same formula to $D^{'}$
to deduce $D'= 0$ since $\delta(D')= 0$. In conclusion, $D= \delta(\sigma_{D})$, and $D$ is exact.
\end{proof}

In the case of transitive Lie
algebroids (i.e. with surjective anchor), we deduce 

\begin{corollary} For any transitive algebroid  $A$,
\[ \de{H}^{*}(A)\cong H^{*}(A; \mathfrak{g}(A)) \ .\]
\end{corollary}

\subsection{The action on De Rham cohomology} 
\label{act-Lie-coh}
As we have already seen in subsection \ref{alternative}, the deformation complex $\de{C}^{*}(A)$ acts on the complex $C^*(A)$ computing the De Rham cohomology of the Lie algebroid $A$. 
More generally, it acts on all complexes $C^{*}(A; E)$ with $E$ a representation of $A$: one keeps the same formulas as in \ref{alternative}, except 
that one replaces the Lie derivatives $L_{\cdot}$ (of smooth functions) in (\ref{nunu}) by the covariant derivatives $\nabla_{\cdot}$ (of sections of $E$). Again, one obtains an action of $\de{C}^{*}(A)$ on $C^{*}(A; E)$. Passing to cohomology, we obtain:

\begin{proposition} 
For any representation $E$ of the Lie algebroid $A$, there is an induced action
\[ \de{H}^{p+1}(A) \otimes H^{q}(A; E)\rmap H^{p+q}(A; E), (D, c)\mapsto L_{D}(c) \]
which makes $H^{*}(A; E)$ into a graded module over the graded Lie algebra $\de{H}^{*}(A)$. 
\end{proposition}

When $D\in \de{H}^{1}(A)$, $L_{D}$ is given by the standard formula (\ref{Lie-derivative}). There is another interesting particular case of this action, which one obtains by looking at deformations of $A$.
In general, any deformation 
$(A_t)$ of $A$ 
induces a ``variation map'' in cohomology
\[ \partial: H^{*}(A)\rmap H^{*+1}(A) \]
as follows: given a cocycle $c\in C^{k}(A)$, one deforms it to a family of cochains $c_t\in C^{*}(A_t)$, one remarks that
\[ \frac{d}{dt}|_{t=0} \delta_{A_{t}}(c_t)\in C^{k+1}(A) \]
is a cocycle (take the derivative at $t= 0$ of $\delta_{A_t}^2(c_t)= 0$), and denote by $\partial([c])$ the
resulting cohomology class. It is not difficult to see that $\partial$ is nothing but the action $L_{c_0}$ on cohomology,
where $c_0\in \de{H}^{2}(A)$ is the cohomology
class induced by the deformation $A_{t}$ of $A$, see Proposition \ref{prop:deformation}. In particular, this shows that the description above
for $\partial$ does not depend on the choices one makes, and it only depends on the equivalence class of the deformation.

\subsection{Poisson manifolds II}
\label{lin-cat}
There is yet another relation between Lie algebroids and Poisson manifolds: Lie algebroid structures on the vector bundle $A$ are in 1-1 correspondence with Poisson structures on its dual $A^{\vee}$ which are ``linear on the fibers'' \cite{Co}. The deformation complex gives a conceptual interpretation (and a new proof) of this result, which further implies that the deformation cohomology of $A$ coincides with the ``lin-Poisson cohomology'' of $A^{\vee}$. Some of the ideas presented in this subsection are closely related to Mackenzie's work on double structures and the work of Mackenzie and Xu on multiplicative vector fields; see \cite{MX} and the references therein). 

First of all, it is useful to change the language slightly, and look at vector bundles as being manifolds with a certain partial linear structure. 
Accordingly, vector bundles will be called $lin$-manifolds, and the category of vector bundles (with varying base!) will be called the $lin$-category.
Hence a $lin$-manifold $E$ has an underlying manifold $E_0$ over which it is a vector bundle. 
Many of the classical objects have a 
corresponding $lin$-version. For instance, given a $lin$-manifold $E$, we can talk about $lin$-vector bundles $\mathcal{E}$ over $E$:
the addition will be a map $\mathcal{E}\times_{E}\mathcal{E}\rmap \mathcal{E}$ in the $lin$ category.  In Mackenzie's terminology, $\mathcal{E}$ will be a double vector bundle: over $E$, and over $\mathcal{E}_{0}$, which are both vector
bundles over $M= E_0$.
The space $\Gamma_{lin}(\mathcal{E})$ is defined as the space of sections $s: E\rmap \mathcal{E}$ which are morphisms 
in the $lin$-category (hence lie over a section $s_0: E_0\rmap \mathcal{E}_0$ of the vector bundle $\mathcal{E}_0$).

There are two important examples: the tangent bundle $TE$ of a $lin$-manifold $E$ is naturally a $lin$-vector bundle over $E$ with 
$(TE)_0= T(E_0)$, while the cotangent bundle $T^{\vee}E$ is a $lin$-vector bundle over $E$
with $(T^{\vee}E)_0= E^{\vee}$ (where the projection
$T^{\vee}E\rmap E^{\vee}$ comes from the inclusion $\pi^*E\subset TE$, $\pi: E\rmap M$ being the projection).
Sections of $TE$ and $T^{\vee}E$ in the $lin$-category will define the space
$\mathcal{X}_{lin}(E)$ of linear vector fields on $E$, and $\Omega^{1}_{lin}(E)$ of linear 1-forms on $E$ (see also \cite{MX}). 
For instance, an element in $\mathcal{X}_{lin}(E)$ consists of a vector field $X$ on $E$, and a vector field $X_{0}$ on $M$, such that 
$X: E\rmap TE$ is a vector bundle map over $X_{0}: M\rmap TM$. Of course, $\mathcal{X}_{lin}(E)\subset \mathcal{X}(E)$, and, locally,
(with respect to local coordinates $x_i$ in $M$ and a basis $e_i$ in $E$), the linear vector fields
are the vector fields on $E$
\[ \sum a_{i} \frac{\partial}{\partial x_i} + \sum b_{j} \frac{\partial}{\partial e_j} \]
with the property that $a_{i}= a_{i}(x)$ depends only on $x\in M$, and $b_{j}= b_{j}(x, v)$ is linear in $v$.
In the same way we can talk about the space $\Gamma_{lin}(\mathcal{E})$ of any $lin$-vector bundle $\mathcal{E}$ over 
$E$, so that 
\[ \mathcal{X}_{lin}(E)= \Gamma_{lin}(TE), \Omega^{1}_{lin}(E)= \Gamma_{lin}(T^{\vee}E) .\]
Note also that $\Omega^{1}_{lin}(E)\subset \Omega^{1}(E)$ consists of those 1-forms $\omega$ on $E$ with the property that
$\omega(X)\in C^{\infty}_{lin}(E)$ for all $X\in \mathcal{X}_{lin}(E)$. We define $\Omega^{k}_{lin}(E)\subset \Omega^{k}(E)$ by the
similar property: when applied to linear vector fields, it must produce linear smooth functions.
And, dually, pairing multivector fields with wedge products of 1-forms, we define the spaces $\mathcal{X}^{k}_{lin}(E)$
of linear $k$-multi-vector fields on $E$. In local coordinates, these are vector fields which are sums of vectors of type
\begin{equation}
\label{local-form} 
a(x)\frac{\partial}{\partial x_{i_1}}\ldots \frac{\partial}{\partial x_{i_{k}}} + \sum b(x, v) \frac{\partial}{\partial x_{i_1}}\ldots \frac{\partial}{\partial x_{i_{k-1}}}\frac{\partial}{\partial e_j},
\end{equation}
with $b(x, v)$ linear in $v$. Note that these spaces are closed under the Nijenhuis-Schouten
bracket on multi-vector fields, denoted by $[\cdot, \cdot]$. 

According to the general philosophy, a linear Poisson structure on a vector bundle $E$ over $M$  
is a Poisson structure whose Poisson tensor $\pi$ is linear:
$\pi \in\mathcal{X}^{2}_{lin}(E)$. Also, the space of linear-multivector fields on $E$, $\mathcal{X}_{lin}(E)$, define
a sub-complex of the Poisson complex of $E$, and we define the linear Poisson cohomology of $E$, denoted $H^{*}_{\pi, lin}(E)$, as
the cohomology of the resulting complex.

Next, any $X\in \mathcal{X}(E^{\vee})$ induces a multi-derivation $D_{X}$ with
\[ D_{X}(s_1, \ldots , s_k)= (ds_1\wedge \ldots \wedge ds_k)(X),\]
where, for $s\in \Gamma(E)$, one views $s$ as a linear smooth function on $E^{\vee}$, and $ds\in \Omega^{1}_{lin}(E^{\vee})$ is its differential.
The symbol $\sigma_X$ of $D_X$ is determined by
\[ d\phi(\sigma_X(s_1, \ldots , s_{k-1}))= (ds_1\wedge \ldots \wedge ds_{k-1}\wedge d\phi)(X),\]
where, as before, for $\phi\in C^{\infty}(M)$, we denote by the same letter the function induced on $E^{\vee}$ (constant on the fibers).
From the local form (\ref{local-form}) of linear multi-vector fields, we see that these expressions determine $X$ completely.
After a careful computation (needed to identify the brackets), one concludes that this construction 
defines an isomorphism of graded Lie algebras (note also that, since the map $X\mapsto D_X$ is local over $M$, one can work in local coordinates)
\[ Der^{k-1}(E)\cong \mathcal{X}^{k}_{lin}(E^{\vee}) \]
(This is a generalization to multivector fields of Proposition 2.4 of \cite{MX}; see also \cite{Yv}).
In particular, we deduce the following proposition, the first part of which is well-known (it first appeared in \cite{Co}, Theorem 2.1.4):

\begin{proposition} Given a vector bundle $A$ over $M$, there is a 1-1 correspondence between Lie algebroid structures on $A$ and 
linear Poisson structures $\pi$ on $A^{\vee}$. Moreover, for any Lie algebroid $A$,
\[ \de{H}^{*}(A)\cong H^{*}_{\pi, lin}(A^{\vee}) .\]
\end{proposition}

\subsection{Rigidity} Finally, we conjecture the following (cohomological) rigidity result similar to known rigidity results for compact Lie groups. The analogue for groupoids of the compactness property of groups is known as properness. For those aspects of Lie groupoids and properness which are relevant to our discussion, we refer to \cite{Cra}.

\begin{conjecture}
If $A$ is a Lie algebroid which admits a proper integrating groupoid $\mathcal{G}$ whose $s$-fibers are $2$-simply connected,
then $\de{H}^{2}(A)= 0$.
\end{conjecture}

Such a result is relevant to the study of smooth deformations of Lie algebroids, and to linearization problems. 

\begin{proposition} 
The conjecture is true for regular algebroids, and for action algebroids. 
\end{proposition}

\begin{proof}
First notice that by Theorem 4 and Proposition 1 of \cite{Cra}, the De Rham cohomology $H^*(A; E)$ vanishes in degrees $1$ and $2$, for any representation $E$. Now the regular case follows by Theorem \ref{regular}, and the action case by 
Proposition \ref{lie-acti}.
\end{proof}

We believe that the previous conjecture can be proven by further working in the $lin$-category (subsection \ref{lin-cat}): one can talk about $lin$-algebroids (the LA-groupoids of Mackenzie \cite{Ma2}), $lin$-groupoids (the VB groupoids of Mackenzie \cite{Ma2}), etc, so that $\de{H}^{*}(A)$ is isomorphic to the ``$lin$-algebroid cohomology'' of the $lin$-algebroid $T^{\vee}A^{\vee}$ (the algebroid associated to the $lin$-Poisson manifold $A^{\vee}$), and then, try to prove a $lin$-version of the van Est-type results of \cite{Cra} (for the $lin$-groupoid $T^{\vee}\mathcal{G}$ over $A^{\vee}$). 



\bibliographystyle{amsplain}
\def\lllll{}

\end{document}